\documentstyle[12pt]{article}
\catcode`\@=11
\@addtoreset{equation}{section}

\catcode`\@=12
\newtheorem{Theorem}{Theorem}[section]
\newtheorem{Definition}{Definition}[section]
\newtheorem{Proposition}{Proposition}[section]
\newtheorem{Lemma}{Lemma}[section]
\newtheorem{Corollary}{Corollary}[section]

\title{A Note On Weinstein Conjecture\thanks{Project 19871044 Supported by NSF}}

\author{Renyi Ma\\
Department of Mathmatics \\
Tsinghua University \\
Beijing, 100084\\
People's Republic of China}

\date { }

\begin{document}
\textwidth=125mm
\textheight=185mm
\parindent=8mm
\frenchspacing
\maketitle

\begin{abstract}
In this article, we give new proofs on  
the some cases on Weinstein conjecture and get some new results on Weinstein conjecture. 
\end{abstract}


\noindent{\bf 2000MR Subject Classification} 32Q65,53D35,53D12

\section{Introduction and results}

Let $\Sigma$ be a smooth closed oriented manifold of dimension
$2n-1$. A contact form on $\Sigma$ is a $1-$form such that
$\lambda \wedge (d\lambda )^{n-1}$ is a volume form on $\Sigma$.
Associated to $\lambda$ there is the so-called Reed vectorfield $X_\lambda $ defined
by $i_{X_\lambda }\lambda   \equiv 1, \ \ i_{X_\lambda }d\lambda  \equiv 0$.
The integral curve of $X_\lambda $ is called $characteritcs$. 
There is a well-known conjecture raised by Weinstein in \cite{we} 
which concerned the close Reeb orbit in 
a contact manifold.

{\bf Conjecture}(see\cite{we}). 
If $(\Sigma ,\lambda )$ is a close simply connected contact manifold 
with contact form $\lambda $ of dimension $2n-1$, then 
there is a close characteristics.

Let $(M,\omega )$ be a symplectic manifold and 
$h(t,x)(=h_t(x))$ a compactly supported smooth function on 
$M\times [0,1]$. Assume that the segment 
$[0,1]$ is endowed with time coordinate $t$. For every function 
$h$ define the $(time-dependent)$ $Hamiltonian$ 
$vector$ $field$ $X_{h_t}$ by the equation:
\begin{equation}
dh_t(\eta )=\omega (\eta ,X_{h_t}) \ \ for \ every \ \eta \in TM
\end{equation}
The flow $g^t_h$ generated by the field 
$X_{h_t}$ is called $Hamiltonian $ $flow$ and its time one map $g_h^1$ is called $Hamiltonian $ $diffeomorphism $.
Now assume that $H$ be a time independent smooth function on 
$M$ and $X_H$ its induced vector field.

\begin{Theorem}
Let $(M,\omega )$ be an exact symplectic manifold convex at infinity or with bounded geometry.  
Let $(\Sigma ,\lambda )$ be a contact manifold 
of induced type  in $M$ with 
induced contact form $\lambda $, i.e., 
there exists a vector field $X$ transversal to $\Sigma $ such that 
$L_X\omega =\omega $ 
and $\lambda =i_X\omega $,  $X_{\lambda } $ its Reeb vector field. 
If there exists a 
Hamiltonian diffeomorphism $h$ such that 
$h(\Sigma )\cap \Sigma =\emptyset $,  
then there exists at least one close characteritics
on $\Sigma $
\end{Theorem}

\begin{Corollary}
Let $(M,\omega )$ be an exact symplectic manifold which is  
convex at infinity or has bounded geometry(see\cite{gro}). 
$M\times C$ be a symplectic manifold with 
symplectic form $\omega \oplus \sigma $, 
here $(C,\sigma )$ standard 
symplectic plane. Let 
$r_0>0$ be a fixed number and 
$B_{r_0}(0)\subset C$ the closed ball with radius 
$r_0$. If 
$(\Sigma ,\lambda )$ be a contact manifold 
of induced type  in $M\times B_{r_0}(0)$ with 
induced contact form $\lambda $, i.e., 
there exists a vector field $X$ transversal to $\Sigma $ such that 
$L_X(\omega \oplus \sigma )=\omega \oplus \sigma $ 
and $\lambda =i_X(\omega \oplus \sigma )$,  $X_{\lambda } $ its Reeb vector 
field. Then 
there exists at least one close characteristics. 
\end{Corollary} 
Corollary 1.1was proved in \cite{ma} by Hofer-Viterb's method(see\cite{hv}).
 
\begin{Corollary}
Let $M$ be any open manifold and $(T^*M,d\alpha )$ be its cotangent bundle.  
Let $(\Sigma ,\lambda )$ be a close contact manifold 
of induced type  in $T^*M$. 
there exists at least one close characteristics on $\Sigma $. 
\end{Corollary}
Corollary 1.2 generalizes the results in \cite{hv1,vi2,ma}.
The proof  
of Theorem1.1 is close as in \cite{mo}.

\section{Lagrangian Non-squeezing}

Let $W$ be a Lagrangian submanifold in $M$, i.e., 
$\omega |W=0$. 

\begin{Definition}
Let 
$$l(M,W, \omega )=\inf \{ |\int _{D^2}f^*\omega |>0 |f:(D^2,\partial D^2)\to (M,W) \} $$
\end{Definition}

\begin{Theorem}
(\cite{po})Let $(M,\omega )$ be a closed  
compact symplectic manifold or a manifold convex at infinity and 
$M\times C$ be a symplectic manifold with 
symplectic form $\omega \oplus \sigma $, 
here $(C,\sigma )$ standard 
symplectic plane. Let 
$2\pi r_0^2<s(M,\omega )$ and 
$B_{r_0}(0)\subset C$ the closed disk with radius 
$r_0$. If $W$ is a close Lagrangian manifold in $M\times B_{r_0}(0)$, then 
$$l(M,W, \omega )<2\pi r_0^2$$
\end{Theorem}
This can be considered as an Lagrangian version of Gromov's symplectic squeezing.

\begin{Corollary}
(Gromov\cite{gro})Let $(V',\omega ')$ be an exact symplectic manifold with 
restricted contact boundary and $\omega '=d\alpha ' $. Let 
$V'\times C$ be a symplectic manifold with 
symplectic form $\omega '\oplus \sigma =d\alpha 
=d(\alpha '\oplus \alpha _0$, 
here $(C,\sigma )$ standard 
symplectic plane. 
If $W$ is a close exact Lagrangian submanifold, then 
$l(V'\times C,W, \omega )==\infty $, i.e., there does not exist any 
close exact Lagrangian submanifold in $V'\times C$.
\end{Corollary}

\begin{Corollary}
Let $L^n$ be a close Lagrangian in $R^{2n}$ and $L(R^{2n},L^n,\omega )=2\pi r_0^2>0$, then 
$L^n$ can not be embedded in $B_{r_0}(0)$ as a Lagrangian submanifold. 
\end{Corollary}

\section{Constructions of Lagrangian submanifolds}

Let $(\Sigma ,\lambda )$ be a contact manifolds with contact form 
$\lambda $ and $X$ its Reeb vector field, then 
$X$ integrates to a Reeb flow $\eta _t$ for $t\in R^1$. 
Let 
$$(V',\omega ')=((R\times \Sigma )\times (R\times \Sigma ), 
d(e^a\lambda )\ominus d(e^b\lambda ))$$ 
and 
$${\cal {L}}=\{ ((0,\sigma ),(0,\sigma ))|(0,\sigma )\in R\times \Sigma \}.$$
Let 
$$L'={\cal {L}}\times R,
L_s'={\cal {L}}\times \{s\}.$$
Then 
define 
\begin{eqnarray}
&&G':L'\to V'\cr 
&&G'(l')=G'(((\sigma ,0),(\sigma ,0)),s)
=((0,\sigma ),(0,\eta _s(\sigma )))
\end{eqnarray}
Then
    
$$W'=G'(L')=\{ ((0,\sigma ),(0,\eta _s(\sigma )))
|(0,\sigma )\in R\times \Sigma ,s\in R\}$$

$$W_s'=G'(L'_s)=\{ ((0,\sigma ),(0,\eta _s(\sigma )))
|(0,\sigma )\in R\times \Sigma \}$$
for fixed $s\in R$.

\begin{Lemma}
There does not exist any Reeb closed orbit in 
$(\Sigma ,\lambda )$ if and only if  
$W'_s\cap W'_{s'}$ is empty for $s\ne s'$.
\end{Lemma}
Proof. First if there exists a closed Reeb orbit in 
$(\Sigma ,\lambda )$, i.e., there exists 
$\sigma _0\in \Sigma $, $t_0>0$ such that 
$\sigma _0=\eta _{t_0}(\sigma _0)$, then 
$((0,\sigma _0),(0,\sigma _0))\in W'_0\cap W'_{t_0}$.  
Second if there exists $s_0\ne s_0'$ such 
that $W'_{s_0}\cap W'_{s_0'}\ne \emptyset $, i.e., 
there exists $\sigma _0$ such that 
$$((0,\sigma _0),(0,\eta _{s_0}(\sigma _0))
=((0,\sigma _0),(0,\eta _{s_0'}(\sigma _0)),$$ 
then $\eta _{(s_0-s_0')}(\sigma _0)=\sigma _0$, i.e., 
$\eta _t(\sigma _0)$ is a closed Reeb orbit.

\begin{Lemma}
If there does not exist any closed Reeb orbit in 
$(\Sigma ,\lambda )$ then 
there exists a smooth Lagrangian injective immersion 
$G':W'\to V'$ with $G'(((0,\sigma ),(0,\sigma )),s)
=((0,\sigma ),(0,\eta _s(\sigma )))$
such that
\begin{equation}
G'_{s_1,s_2}:{\cal L}\times (-s_1,s_2)\to V'
\end{equation}
is a regular exact Lagrangian embedding for any finite real number 
$s_1$, $s_2$, here we denote by $W'(s_1,s_2)=G'_{s_1,s_2}({\cal L}\times 
(s_1,s_2))$. 
\end{Lemma}
Proof. One check 
\begin{equation}
{G'}^*((e^a\lambda -e^b\lambda ))
=\lambda -\eta (\cdot ,\cdot )^*\lambda 
=\lambda -(\eta _s^*\lambda +i_X\lambda ds)=-ds
\end{equation}
since $\eta _s^*\lambda =\lambda $. 
This implies that ${G}'$ is an exact  Lagrangian embedding, this proves 
Lemma 3.2. 

\vskip 3pt 

Now we modify the above construction as follows:  
\begin{eqnarray}
&&F':{\cal {L}}\times R\times R\to (R\times \Sigma )\times 
(R\times \Sigma )\cr 
&&F'(((0,\sigma ),(0,\sigma )),s,b)=((0,\sigma ),(b,\eta _s(\sigma )))
\end{eqnarray}
Now we embed a elliptic curve $E$ long along $s-axis$ and thin along $b-axis$ such that 
$E\subset [-s_1,s_2]\times [0,\varepsilon]$. We parametrize the $E$ by $t$.

\begin{Lemma}
If there does not exist any closed Reeb orbit 
in $(\Sigma ,\lambda )$, 
then 
\begin{eqnarray}
&&F:{\cal {L}}\times S^1\to (R\times \Sigma )\times 
(R\times \Sigma )\cr 
&&F(((0,\sigma ),(0,\sigma )),t)=((0,\sigma ),(b(t),\eta _{s(t)}(\sigma )))
\end{eqnarray}
is a compact Lagrangian submanifold. Moreover 
\begin{equation}
l(V',F({{\cal {L}}}\times S^1,d(e^a\lambda -e^b\lambda ))=area(E)
\end{equation}
\end{Lemma}
Proof. We check
that 
\begin{eqnarray}
{F}^*(e^a\lambda \ominus e^b\lambda )&=&-e^{b(t)}ds(t)
\end{eqnarray}
So, $F$ is a Lagrangian embedding.

If the circle $C$ homotopic to $C_1\subset {\cal {L}}\times s_0$ then  we compute
\begin{eqnarray}
\int _CF^*(e^a\lambda -e^b\lambda )=\int _{C_1}F^*(e^0\lambda -e^0\lambda )=0. 
\end{eqnarray}
since $\lambda-\lambda  |C_1=0$ due to $C_1\subset {\cal {L}}$.
If the circle $C$ homotopic to $C_1\subset l_0\times S^1$ then  we compute
\begin{eqnarray}
\int _CF^*(e^a\lambda -e^b\lambda )=\int _{C_1}(-)e^bds=n(area(E)). 
\end{eqnarray}
This proves the Lemma.

{\bf Gromov's figure eight construction:}
First we note that the construction of section 3.1 holds for any symplectic manifold.
Now let $(M,\omega )$ be an exact symplectic manifold with 
$\omega =d\alpha $. 
Let $\Sigma =H^{-1}(0)$ be a regular and close smooth 
hypersurface in $M$. $H$ is a time-independent Hamilton function. 
Set 
$(V',\omega ')=(M\times M,\omega \ominus \omega )$.
If there does not exist any close 
orbit for $X_H$ 
in $(\Sigma ,X_H)$, one can construct 
the Lagrangian 
submanifold 
$L$ as in section 3.1, let $W'=L$.
Let $h_t=h(t,\cdot ):M\to M$, $0\leq t\leq 1$ be a Hamiltonian isotopy 
of $M$ induced by hamilton fuction $H_t$ such that $h_1(\Sigma )\cap \Sigma =\emptyset $, 
$|H_t|\leq C_0$.
Let $\bar h_t=(id,h_t)$. Then $F'_t=\bar h_t:W'\to V'$ 
be an isotopy of Lagrangian embeddings. 
As in \cite{gro}, we can use symplectic figure eight trick invented by Gromov to 
construct a Lagrangian submanifold $W$ in $V=V'\times R^2$ through the 
Lagrange isotopy $F'$ in $V'$, i.e., we have 
\begin{Proposition}
Let $V'$, $W'$ and $F'$ as above. 
Then there exists  
a weakly exact Lagrangian embedding $F:W'\times S^1\to V'\times R^2$ 
with $W=F(W'\times S^1)$ is contained in 
$M\times M\times B_R(0)$, here $4\pi R^2=8C_0$ and 
\begin{equation}
l(V',W,\omega )=area(M'_0)=A(T).
\end{equation}
\end{Proposition}
Proof. Similar to \cite[2.3$B_3'$]{gro}.

{\bf Example.} Let $M$ be an open manifold and $(T^*M,p_idq_i)$ be 
the cotangent bundle of open manifold with the 
Liouville form $p_idq_i$. Since 
$M$ is open, there exists a function $g:M\to R$ without 
critical point. The translation by 
$tTdg$ along the fibre gives a hamilton isotopy 
of $T^*M:h^T_t(q,p)=(q,p+tTdg(q))$, so
for any given compact set $K\subset T^*M$, there exists 
$T=T_K$ such that $h^T_1(K)\cap K=\emptyset $.

\subsection{Proof on Theorem 1.1}

Since $(\Sigma ,\lambda )$ be a close contact manifold 
of induced type  in $M$ with 
induced contact form $\lambda $, then by the well known theorem that 
the neighbourhood $(U(\Sigma ),\omega )$ of $\Sigma $ is symplectomorphic to 
$([-\varepsilon ,\varepsilon]\times \Sigma ,de^a\lambda )$ for small 
$\varepsilon $. So, by Proposition 3.1, we have 
a close Lagrangian submanifold $F({\cal {L}}\times S^1)$ contained in 
$M\times M\times B_R(0)$. By Lagrangian squeezing theorem, i.e., Theorem 2.1, 
we have 
\begin{equation}
l((M\times M\times C),F({{\cal {L}}}\times S^1,\omega \oplus \omega )=area(E)\leq 2\pi R^2.
\end{equation}
If $s_2-s_1$ large enough, $area(E)>2\pi R^2$. This 
is a contradiction. This contradiction shows 
there exists at least one close characteristics.

\end{document}